\documentclass[12pt,twoside]{article}
\usepackage{amsmath,amsbsy,amsfonts}

\newtheorem{lemma}{Lemma}[section]

\newtheorem{theorem}{Theorem}[section]

\newtheorem{definition}{Definition}[section]

\let\Section=\section
\def\section{\setcounter{equation}{0}\Section}
\input{tcilatex}
\begin{document}

\title{\textbf{Boundedness and blow-up of solutions for a nonlinear elliptic
system}}
\author{\vspace{1mm}Dragos-Patru Covei \\
\textsf{\ }{\small \textit{Department of Development}}\\
{\small \textit{Constantin Brancusi University of Tg-Jiu, Gorj, Romania}}\\
}
\maketitle

\begin{abstract}
The main objective in the present paper is to obtain the existence results
for bounded and unbounded solutions of some quasilinear elliptic systems.
Related results as obtained here have been established recently in [C. O.
Alves and A. R.F. de Holanda, \textit{Existence of blow-up solutions for a
class of elliptic systems}, \ Differential Integral Equations, Volume 26,
Number 1/2 (2013), Pages 105-118.]. Also, we present some references to give
the connection between these type of problems with probability and
stochastic processes, hoping that these are interesting for the audience of
analysts likely to read this paper.
\end{abstract}

\section{Introduction}

The question of existence of solutions for elliptic equation of the form%
\begin{equation}
\Delta _{p}u=f\left( x,u\right) \text{ in }\Omega ,  \label{ko}
\end{equation}%
was studied by many researchers (see Bandle and Marcus \cite{bandle}, the
author \cite{covei}, Lair \cite{lair}, Matero \cite{matero}, Mohammed \cite%
{ahmed}, Peterson-Wood \cite{w} with their references). This work is devoted
to the study of the more general nonlinear elliptic problems of the type%
\begin{equation}
\left\{ 
\begin{array}{l}
\Delta _{p}u_{i}=F_{u_{i}}\left( x,u_{1},...,u_{i},...,u_{d}\right) \text{
in }\Omega , \\ 
i=1,...,d%
\end{array}%
\right.  \label{P}
\end{equation}%
where $d\geq 1$ is integer, $\Omega \subset \mathbb{R}^{N}\left( N>1\right) $
is a smooth, bounded domain or $\Omega =\mathbb{R}^{N}$ , $\Delta _{p}u_{i}:=%
\func{div}\left( \left\vert \nabla u_{i}\right\vert ^{p-2}\nabla
u_{i}\right) $ ($1<p<\infty $) is the p-Laplacian operator and $F_{u_{i}}$ ($%
i=1,...,d$) stands for the derivatives of a continuously differentiable
function $F:\Omega \times \lbrack \mathbb{R}^{+}]^{d}\rightarrow \mathbb{R}%
^{+}$ in $\left( u_{1},...,u_{d}\right) $. For the case $\Omega =\mathbb{R}%
^{N}$, we also consider the following class of elliptic systems:%
\begin{equation}
\left\{ 
\begin{array}{l}
\Delta _{p}u_{i}=a_{i}\left( x\right) F_{u_{i}}\left(
x,u_{1},...,u_{i},...,u_{d}\right) \text{ in }\mathbb{R}^{N}, \\ 
u_{i}>0\text{ in }\mathbb{R}^{N}, \\ 
i=1,...,d%
\end{array}%
\right.  \label{LS}
\end{equation}%
where $a_{i}:\mathbb{R}^{N}\rightarrow \left( 0,\infty \right) $ are
suitable functions. Associated with the class of systems (\ref{LS}), our
main result is concerned with the existence of entire large solutions, that
is, solutions $\left( u_{1},...,u_{d}\right) $ satisfying $u_{i}\left(
x\right) \rightarrow \infty $ as $\left\vert x\right\vert \rightarrow \infty 
$ for all $i=1,...,d$.

The interest on systems (\ref{P})-(\ref{LS}) comes from some problems
studied in the works of Lasry-Lions \cite{lasry}, Busca-Sirakov \cite{busca}
and Dynkin \cite{dynkin} where the authors give the connection between these
type of problems with probability and stochastic processes and from the
recently work of Alves and Holanda \cite{alves} where these systems are
considered for the case $p=2$, in terms of the pure mathematics. The
difference between our work and the paper by \cite{alves} is that: our
systems can have any number of equations, the potential functions $a_{i}$
cover more general properties and that we use in the proofs theories for
quasilinear operators instead of the theories for linear operators used by
Alves-Holanda \cite{alves}. We also remark that th authors Alves-Holanda 
\cite{alves} extended the results of Bandle-Marcus \cite{bandle}, obtained
for the scalar equation in bounded domains, to the system of two equations
while our proof work for any numbers of equations.

To begin with our results we make the following convention: we say that a
function $h:\left[ 0,+\infty \right) \rightarrow \left[ 0,\infty \right) $
belongs to $\mathcal{F}$ if 
\begin{eqnarray*}
h &\in &C^{1}\left( \left[ 0,\infty \right) \right) ,\text{ }h\left(
0\right) =0,\text{ }h^{\prime }\left( t\right) \geq 0\text{ }\forall t\in %
\left[ 0,\infty \right) , \\
h\left( t\right) &>&0\text{ }\forall t\in \left( 0,\infty \right)
\end{eqnarray*}%
and the Keller-Osserman \cite{keller}, \cite{osserman} condition is
satisfied, that is,%
\begin{equation*}
\int_{1}^{\infty }\frac{1}{H\left( t\right) ^{1/p}}dt<\infty ,
\end{equation*}%
where $H\left( t\right) =\int_{0}^{t}h\left( s\right) ds$.

Our main result for problem (\ref{P}) on a bounded domain is the following:
\ 

\begin{theorem}
\label{1.1}Suppose $\Omega $ is a smooth, bounded domain in $\mathbb{R}^{N}$
and that there exist $f_{i},$ $g\in \mathcal{F}$ satisfying%
\begin{equation}
F_{t_{i}}\left( x,t_{1},...,t_{i},...,t_{d}\right) \geq f_{i}\left(
t_{i}\right) \text{ }\forall x\in \overline{\Omega },\text{ }t_{i}>0\text{
and }i=1,...,d  \label{1.2}
\end{equation}%
and%
\begin{equation}
g\left( t\right) \geq \max \left\{ F_{t}\left( x,t,...,t\right) \right\} 
\text{ }\forall x\in \overline{\Omega },\text{ }t>0.  \label{1.4}
\end{equation}%
Then:

1. problem (\ref{P}) admits a positive solution with boundary condition 
\begin{equation}
\left\{ 
\begin{array}{l}
u_{i}=\alpha _{i}\text{ on }\partial \Omega \\ 
i=1,...,d \\ 
\alpha _{i}\in \left( 0,\infty \right) .%
\end{array}%
\right.  \label{F}
\end{equation}

2. problem (\ref{P}) admits a positive solution with the boundary condition 
\begin{equation}
\left\{ 
\begin{array}{l}
u_{i}=\infty \text{ on }\partial \Omega \\ 
i=1,...,d%
\end{array}%
\right.  \label{I}
\end{equation}%
where $u_{i}=\infty $ on $\partial \Omega $ should be understood as $%
u_{i}\left( x\right) \rightarrow \infty $ as $dis\left( x,\partial \Omega
\right) \rightarrow 0$.

3. problem (\ref{P}) admits a positive solution with boundary condition:
there are $i_{0},j_{0}\in \left\{ 1,...,d\right\} $ such that 
\begin{equation}
\left\{ 
\begin{array}{l}
\text{ }u_{i_{0}}=\infty \text{ on }\partial \Omega \text{ } \\ 
\text{ }u_{j_{0}}<\infty \text{ on }\partial \Omega \text{ for any }%
j_{0}\neq i_{0}\text{ }%
\end{array}%
\right.  \label{SF}
\end{equation}%
and the set $\left\{ 1,...,d\right\} $ is crossed by $i_{0}$ respectively $%
j_{0}$.
\end{theorem}

Our next result is related to the existence of a solution for system (\ref%
{LS}). For expressing the next result, we assume that functions $a_{i}$ ($%
i=1,...,d$) satisfy the following conditions: 
\begin{equation}
a_{i}\left( x\right) >0\text{ for all }x\in \mathbb{R}^{N}\text{ and }%
a_{i}\in C_{loc}^{0,\vartheta }\left( \mathbb{R}^{N}\right) ,\vartheta \in
\left( 0,1\right)  \label{1.5.}
\end{equation}%
and that the quasilinear system%
\begin{equation}
\begin{array}{c}
-\Delta _{p}z\left( x\right) =\overset{d}{\underset{i=1}{\sum }}a_{i}\left(
x\right) \text{ for }x\in \mathbb{R}^{N}\text{, \ }z\left( x\right)
\rightarrow 0\text{ as }\left\vert x\right\vert \rightarrow \infty%
\end{array}
\label{1.6.}
\end{equation}%
has a $C^{1}$-upper solution, in the sense that%
\begin{equation}
\left\{ 
\begin{array}{l}
\int_{\mathbb{R}^{N}}\left\vert \nabla z\right\vert ^{p-2}\nabla z\cdot
\nabla \phi dx\geq \int_{\mathbb{R}^{N}}\sum_{i=1}^{d}a_{i}\left( x\right)
\phi dx\text{, }\phi \in C_{0}^{\infty }\left( \mathbb{R}^{N}\right) \text{, 
}\phi \geq 0 \\ 
z\in C^{1}\left( \mathbb{R}^{N}\right) \text{, }z>0\text{ in }\mathbb{R}^{N}%
\text{, }z\left( x\right) \rightarrow 0\text{ as }\left\vert x\right\vert
\rightarrow \infty .%
\end{array}%
\right.  \label{qa}
\end{equation}

\begin{theorem}
\label{entire}Assume that (\ref{1.2})-(\ref{1.4}), (\ref{1.5.})-(\ref{qa})
hold. Then system (\ref{LS}) has an entire large $C^{1}$-solution (in the
distribution sense).
\end{theorem}

To prepare for proving our theorems, we need some additional results.

\section{Preliminary results}

Let $\Omega \subset \mathbb{R}^{N}$ ($N\geq 2$) be a smooth, bounded domain
in $\mathbb{R}^{N}$ and $1<p<\infty $. The first auxiliary result can be
seen in the paper of Matero \cite[pp. 233]{jerk}.

\begin{lemma}
\label{jm}Assume that $g$ meets the conditions: $g$ is a continuous,
positive, increasing function on $\mathbb{R}_{+}$, and $g\left( 0\right) =0$%
. Let $h\in W^{1,p}\left( \Omega \right) $ be such that $\left( G\circ
h\right) \in L^{1}\left( \Omega \right) $ where $G\left( s\right)
=\int_{0}^{s}g\left( t\right) dt$. Then there exists a unique $u\in
W^{1,p}\left( \Omega \right) $ which (weakly) solves the problem%
\begin{equation*}
\left\{ 
\begin{array}{l}
\Delta _{p}u\left( x\right) =g\left( u\left( x\right) \right) ,\text{ }x\in
\Omega , \\ 
u\left( x\right) =h\left( x\right) ,\text{ }x\in \partial \Omega .%
\end{array}%
\right.
\end{equation*}%
Furthermore, if $u_{1}$ and $u_{2}$ are the solutions corresponding to $%
h_{1} $ and $h_{2}$ with $h_{1}\leq h_{2}$ on $\partial \Omega $, then $%
u_{1}\leq u_{2}$ in $\Omega $. Finally, there exists an $\beta \in \left(
0,1\right) $ such that $u\in C^{1,\beta }\left( D\right) $ for any compact
set $D\subset \Omega $.
\end{lemma}

The following comparison principle is proved in the article of Sakaguchi 
\cite{SS} (or consult some ideas of the proof in work of Tolksdorf \cite[%
Lemma 3.1.]{tolk}).

\begin{lemma}
\label{cm}Let $u,v\in W^{1,p}\left( \Omega \right) $ satisfy $-\Delta
_{p}u\leq -\Delta _{p}v$ for $x\in \Omega $, \textit{in the weak sense. }If $%
u\leq v$ on $\partial \Omega $ then $u\leq v$ in $\Omega $.
\end{lemma}

The following Lemma can be found in Zuodong (\cite{yang}).

\begin{lemma}
\label{bo} Suppose $f\in \mathcal{F}$ and that $u\in W_{loc}^{1,p}\left(
\Omega \right) \cap C\left( \Omega \right) $ satisfies 
\begin{equation*}
-\int_{\Omega }\left\vert \nabla u\right\vert ^{p-2}\nabla u\nabla \varphi
dx=\int_{\Omega }f\left( u\right) \varphi dx,\forall \varphi \in
C_{0}^{\infty }\left( \Omega \right) .
\end{equation*}%
Then, there exists a monotone decreasing function $\mu :\left( 0,\infty
\right) \rightarrow \left( 0,\infty \right) $ determined by $f$ such that 
\begin{equation*}
u\left( x\right) \leq \mu \left( dist\left( x,\partial \Omega \right)
\right) \text{ }\forall x\in \Omega .
\end{equation*}%
Moreover,%
\begin{equation*}
\lim_{t\rightarrow 0}\mu \left( t\right) =\infty ,\lim_{t\rightarrow \infty
}\mu \left( t\right) =-\infty .
\end{equation*}
\end{lemma}

Next, we begin with recalling the definition of sub and super-solution used
in the present context. The system that we will study is the following%
\begin{equation}
\left\{ 
\begin{array}{l}
\Delta _{p}u_{i}=G_{u_{i}}\left( x,u_{1},...,u_{i},...,u_{d}\right) \text{
in }\Omega \\ 
u_{i}=f_{i}\text{ on }\partial \Omega \\ 
i=1,...,d%
\end{array}%
\right.  \label{2.1}
\end{equation}%
where $f_{i}\in W^{1,p}\left( \Omega \right) $ and $G\left(
x,t_{1},...,t_{i},...,t_{d}\right) :\Omega \times \left[ \mathbb{R}\right]
^{d}\rightarrow \mathbb{R}$ is measurable in $x\in \Omega $, continuously
differentiable in $t_{i}\in \mathbb{R}$, and satisfies the following
condition: for each $T_{i}>0$ fixed ($i=1,...,d$), there exists $C=C\left(
T_{i}\right) >0$ such that%
\begin{equation}
\left\vert G\left( x,t_{1},...,t_{d}\right) \right\vert \leq C\text{ }%
\forall \left( x,t_{1},...,t_{d}\right) \in \Omega \times \left[ -T_{i},T_{i}%
\right] ^{d}.  \label{2.2}
\end{equation}%
Now we introduce the concept of sub- and super-solution in the weak sense.

\begin{definition}
By definition $\left( \underline{u}_{1},...,\underline{u}_{d}\right) \in %
\left[ W^{1,p}\left( \Omega \right) \right] ^{d}$ is a (weak) sub-solution
to (\ref{2.1}), if $\underline{u}_{i}\leq f_{i}$ on $\partial \Omega ,$%
\begin{equation*}
\int_{\Omega }\left\vert \nabla \underline{u}_{i}\right\vert ^{p-2}\nabla 
\underline{u}_{i}\nabla \phi dx+\int_{\Omega }G_{\underline{u}_{i}}\left( x,%
\underline{u}_{1},...,\underline{u}_{d}\right) \phi dx\leq 0
\end{equation*}%
for all $\phi \in C_{0}^{\infty }\left( \Omega \right) $ with $\phi \geq 0$
and $i=1,...,d$.

Similarly $\left( \overline{u}_{1},...,\overline{u}_{d}\right) \in \left[
W^{1,p}\left( \Omega \right) \right] ^{d}$ \ is a (weak) super-solution to (%
\ref{2.1}) if in the above the reverse inequalities hold.
\end{definition}

The following result holds:

\begin{lemma}
\label{ss}Suppose $\left( \underline{u}_{1},...,\underline{u}_{d}\right) $
is a sub-solution while $\left( \overline{u}_{1},...,\overline{u}_{d}\right) 
$ is a super-solution to problem (\ref{2.1}), and assume that there are
constants $\underline{a}_{i},\overline{a}_{i}\in \mathbb{R}$ such that 
\begin{equation*}
\underline{a}_{i}\leq \underline{u}_{i}\leq \overline{u}_{i}\leq \overline{a}%
_{i}\text{ almost every where in }\Omega .
\end{equation*}%
If (\ref{2.2}) holds, then there exists a weak solution $\left(
u_{1},...,u_{d}\right) \in \left[ W^{1,p}\left( \Omega \right) \right] ^{d}$
of (\ref{2.1}), satisfying the condition%
\begin{equation*}
\underline{u}_{i}\leq u_{i}\leq \overline{u}_{i}\text{ almost everywhere in }%
\Omega .
\end{equation*}
\end{lemma}

We will not give the proof here since he can now proved as in \cite[Theorem
2.1, pp. 110]{alves} with some ideas from \cite{zhou}.

\section{Proof of main results}

In this section, we will prove the main results of this paper.

\subsection{Proof of Theorem \protect\ref{1.1}}

\subsubsection{Proof of 1:}

In what follows, we denote by $\psi \in W^{1,p}\left( \Omega \right) $ the
unique positive solution of the problem%
\begin{equation*}
\left\{ 
\begin{array}{l}
-\int_{\Omega }\left\vert \nabla \psi \right\vert ^{p-2}\nabla \psi \nabla
\phi dx=\int_{\Omega }g\left( \psi \right) \phi dx\text{ in }\Omega ,\text{ }%
\forall \phi \in C_{0}^{\infty }\left( \Omega \right) \text{ with }\phi \geq
0, \\ 
\psi >0\text{ in }\Omega \\ 
\psi =m\text{ on }\partial \Omega%
\end{array}%
\right.
\end{equation*}%
where $m=\min \left\{ \alpha _{1},...,\alpha _{d}\right\} $, which exists
and minimizes the Euler-Lagrange functional%
\begin{equation*}
J\left( \psi \right) =\int_{\Omega }\left( \frac{1}{p}\left\vert \nabla \psi
\right\vert ^{p}+G\left( \psi \left( x\right) \right) \right) dx
\end{equation*}%
on the set%
\begin{equation*}
K=\left\{ v\in L^{1}\left( \Omega \right) \left\vert v-m\in
W_{0}^{1,p}\left( \Omega \right) \text{ and }\left( G\circ v\right) \in
L^{1}\left( \Omega \right) \right. \right\}
\end{equation*}%
i.e., $\psi $ meets the boundary condition $\left( \psi -m\right) \in
W_{0}^{1,p}\left( \Omega \right) $ in the weak sense [see Lemma \ref{jm} and 
\cite[Paragraph 6, pp. 12]{bandle}]. Then, 
\begin{equation*}
\left\{ 
\begin{array}{l}
-\int_{\Omega }\left\vert \nabla \psi \right\vert ^{p-2}\nabla \psi \nabla
\phi dx=\int_{\Omega }g\left( \psi \right) \phi dx\geq \int_{\Omega }F_{\psi
}\left( x,\psi ,...,\psi \right) \phi dx\text{ in }\Omega , \\ 
\psi =m\leq \alpha _{i}\text{ on }\partial \Omega \\ 
i=1,...,d%
\end{array}%
\right.
\end{equation*}%
and so $\left( \underline{u}_{1},...,\underline{u}_{d}\right) =\left( \psi
,...,\psi \right) $ is a sub-solution for the system%
\begin{equation}
\left\{ 
\begin{array}{l}
\Delta _{p}u_{i}=F_{u_{i}}\left( x,u_{1},...,u_{i},...,u_{d}\right) \text{
in }\Omega \\ 
u_{i}=\alpha _{i}\text{ on }\partial \Omega \\ 
i=1,...,d.%
\end{array}%
\right.  \label{PP}
\end{equation}%
Clearly, $\left( \overline{u}_{1},...,\overline{u}_{d}\right) =\left(
M,...,M\right) ,$ with 
\begin{equation*}
M=\max \left\{ \alpha _{i}\left\vert i=1,...,d\right. \right\} ,
\end{equation*}%
is a super-solution of (\ref{PP}). We prove that, $\underline{u}_{i}\leq 
\overline{u}^{i}$ for all $i=1,...,d$. Indeed, 
\begin{equation*}
\left\{ 
\begin{array}{l}
-\Delta _{p}\underline{u}_{i}=-g\left( \underline{u}_{i}\right) \leq -\Delta
_{p}\overline{u}_{i}=0\text{ in }\Omega , \\ 
\underline{u}_{i}=m\leq \overline{u}_{i}=M\text{ on }\partial \Omega , \\ 
i=1,...,d%
\end{array}%
\right.
\end{equation*}%
and then with the use of Lemma \ref{cm} it follows that $\underline{u}%
_{i}\leq \overline{u}_{i}$ in $\overline{\Omega }$. Then, there exists a
critical point $\left( u_{1},...,u_{d}\right) \in \left[ W^{1,p}\left(
\Omega \right) \right] ^{d}$, provided by Lemma \ref{ss}, which minimize the
Euler-Lagrange functional%
\begin{equation*}
I\left( u_{1},...,u_{d}\right) =\frac{1}{p}\int_{\Omega
}\sum_{i=1}^{d}\left\vert \nabla u_{i}\right\vert ^{p}dx+\int_{\Omega
}F\left( x,u_{1},...,u_{d}\right) dx
\end{equation*}%
and that solve, in the weak sense, the system%
\begin{equation}
\left\{ 
\begin{array}{l}
\Delta _{p}u_{i}=F_{u_{i}}\left( x,u_{1},...,u_{i},...,u_{d}\right) \text{
in }\Omega \\ 
u_{i}=\alpha _{i}\text{ on }\partial \Omega \\ 
i=1,...,d%
\end{array}%
\right.  \tag{$P_{\alpha }$}  \label{n}
\end{equation}%
and satisfying $\psi \leq u_{i}\leq M$ in $\Omega $ for all $i=1,...,d$.
Since $u_{i}\in L_{loc}^{\infty }\left( \Omega \right) $, by the regularity
theory \cite{dib,lib,tolk2}, it follows that $u_{i}\in C^{1}\left( \Omega
\right) $.

\subsubsection{Proof of 2:}

To study this case, we begin considering the system%
\begin{equation}
\left\{ 
\begin{array}{l}
\Delta _{p}u_{i}=F_{u_{i}}\left( x,u_{1},...,u_{i},...,u_{d}\right) \text{
in }\Omega , \\ 
u_{i}=n\text{ on }\partial \Omega , \\ 
i=1,...,d.%
\end{array}%
\right.  \label{in}
\end{equation}%
Then, by the finite case above, problem (\ref{in}) has a solution $\left(
u_{1}^{n},...,u_{d}^{n}\right) $.

We prove that the sequence of solutions $\left(
u_{1}^{n},...,u_{d}^{n}\right) $ can be chosen satisfying the inequality%
\begin{equation}
u_{i}^{n}\leq u_{i}^{n+1}\text{ for all }i=1,...,d\text{ and }n\in \mathbb{N}%
.  \label{3.1}
\end{equation}%
To prove this, we consider the solution $\left(
u_{1}^{1},...,u_{d}^{1}\right) $ of the problem 
\begin{equation}
\left\{ 
\begin{array}{l}
\Delta _{p}u_{i}=F_{u_{i}}\left( x,u_{1},...,u_{i},...,u_{d}\right) \text{
in }\Omega , \\ 
u_{i}=1\text{ on }\partial \Omega , \\ 
i=1,...,d.%
\end{array}%
\right.  \label{i1}
\end{equation}%
and note that it is a sub-solution of%
\begin{equation}
\left\{ 
\begin{array}{l}
\Delta _{p}u_{i}=F_{u_{i}}\left( x,u_{1},...,u_{i},...,u_{d}\right) \text{
in }\Omega , \\ 
u_{i}=2\text{ on }\partial \Omega , \\ 
i=1,...,d.%
\end{array}%
\right.  \label{i2}
\end{equation}%
while the pair $(M_{1},...,M_{1})$ is a super-solution of (\ref{i2}) for $%
M_{1}=2$. Once $0\leq u_{i}\left( x\right) \leq 2$ ($i=1,...,d$) $\forall
x\in \overline{\Omega }$ , Lemma \ref{ss} implies that there exists a
solution $\left( u_{1}^{2},...,u_{d}^{2}\right) $ of%
\begin{equation*}
\left\{ 
\begin{array}{l}
\Delta _{p}u_{i}=F_{u_{i}}\left( x,u_{1},...,u_{i},...,u_{d}\right) \text{
in }\Omega , \\ 
u_{i}=2\text{ on }\partial \Omega , \\ 
i=1,...,d.%
\end{array}%
\right.
\end{equation*}%
satisfying $u_{i}^{1}\left( x\right) \leq u_{i}^{2}\left( x\right) $. Using
the argument above, for each $M_{n}=n+1$; $n=1,2,...$, we get a solution $%
\left( u_{1}^{n},...,u_{d}^{n}\right) $ of (\ref{in}), which is a
sub-solution, and the pair $(M_{n},...,M_{n})$ is a super-solution
respectively of 
\begin{equation*}
\left\{ 
\begin{array}{l}
\Delta _{p}u_{i}=F_{u_{i}}\left( x,u_{1},...,u_{i},...,u_{d}\right) \text{
in }\Omega , \\ 
u_{i}=n+1\text{ on }\partial \Omega ,, \\ 
i=1,...,d.%
\end{array}%
\right.
\end{equation*}%
Thereby, the sequence of solutions $\left( u_{1}^{n},...,u_{d}^{n}\right) $
satisfies the inequality (\ref{3.1}). Finally, we construct an upper bound
of the sequence. More exactly, we show that $\left\{ \left(
u_{1}^{n},...,u_{d}^{n}\right) \right\} _{n\geq 1}$ is uniformly bounded in
any compact subset of $\Omega $. To this end, we begin recalling that by (%
\ref{1.2})%
\begin{equation*}
\left\{ 
\begin{array}{l}
\Delta _{p}u_{i}^{n}\geq f_{i}\left( u_{i}^{n}\right) \text{ in }\Omega \\ 
u_{i}^{n}>0\text{ in }\Omega \\ 
u_{i}^{n}\leq n\text{ on }\partial \Omega%
\end{array}%
\right.
\end{equation*}%
with $f_{i}\in \mathcal{F}$. If $\widetilde{u}_{i}^{n}$ ($i=1,...,d$) denote
the unique solutions of the problems%
\begin{equation*}
\left\{ 
\begin{array}{l}
\Delta _{p}u_{i}=f_{i}\left( u_{i}\right) \text{ in }\Omega \\ 
u_{i}>0\text{ in }\Omega \\ 
u_{i}=n\text{ on }\partial \Omega%
\end{array}%
\right.
\end{equation*}%
it follows from Lemma \ref{cm} that%
\begin{equation*}
u_{i}^{n}\leq \widetilde{u}_{i}^{n}\text{ in }\Omega \text{ for all }n\geq 1.
\end{equation*}%
By Lemma \ref{bo}, there exist non-increasing continuous functions $\mu _{i}:%
\mathbb{R}^{+}\rightarrow \mathbb{R}^{+}$ such that%
\begin{equation*}
\widetilde{u}_{i}^{n}\text{ }\leq \mu _{i}\left( dist\left( x,\partial
\Omega \right) \right) \text{ }\forall n\in \mathbb{N},\forall x\in \Omega 
\text{ and }i=1,...,d
\end{equation*}%
showing that%
\begin{equation}
0<u_{i}^{1}\left( x\right) \leq u_{i}^{n}\left( x\right) \leq \mu _{i}\left(
d\left( x\right) \right) \text{ }\forall n\in \mathbb{N},\forall x\in \Omega
\label{3.2}
\end{equation}%
where $d\left( x\right) =dist\left( x,\partial \Omega \right) $. Thus there
exists a subsequence, still denoted again by $u_{i}^{n}$, which converges to
a function $u_{i}$ in $W^{1,p}\left( \Omega \right) $. In other words%
\begin{equation*}
u_{i}\left( x\right) :=\lim_{n\rightarrow \infty }u_{i}^{n}\left( x\right) 
\text{ for all }x\in \Omega \text{ and }i=1,...,d\text{.}
\end{equation*}%
The estimates (\ref{3.2}) combined with the bootstrap argument yield that $%
u_{i}^{n}\left( x\right) \rightarrow u_{i}\left( x\right) $ in $C^{1}\left(
K\right) $ for any compact subset $K\subset \Omega $. \ Furthermore, it is
clear that, $u_{i}\left( x\right) \in C^{1}\left( \Omega \right) $ and $%
\left( u_{1},...,u_{d}\right) $ is a solution of (\ref{P}); that is,%
\begin{equation*}
\left\{ 
\begin{array}{l}
\Delta _{p}u_{i}=F_{u_{i}}\left( x,u_{1},...,u_{i},...,u_{d}\right) \text{
in }\Omega , \\ 
u_{i}>0\text{ in }\Omega , \\ 
i=1,...,d.%
\end{array}%
\right.
\end{equation*}%
To complete the proof, it suffices to prove that $\left(
u_{1},...,u_{d}\right) $ blows up at the boundary. Supposing for the sake of
contradiction that $u_{i}$ does not blow up at the boundary, there exist $%
x_{0}\in \partial \Omega $ and $\left( x_{k}\right) \subset \Omega $\ such
that%
\begin{equation*}
\lim_{k\rightarrow \infty }x_{k}=x_{0}\text{ and }\lim_{k\rightarrow \infty
}u_{i}\left( x_{k}\right) =L_{i}\in \left( 0,\infty \right) .
\end{equation*}%
In what follows, fix $n>4L_{i}$ and $\delta >0$ such that $u_{i}^{n}\left(
x\right) \geq n/2$ for all $x\in \overline{\Omega }_{\delta },$ where 
\begin{equation*}
\overline{\Omega }_{\delta }=\left\{ x\in \overline{\Omega }\left\vert
dist\left( x,\partial \Omega \right) \leq \delta \right. \right\} .
\end{equation*}%
Then, for $k$ large enough, $x_{k}\in \overline{\Omega }_{\delta }$ and $%
u_{i}^{n}\left( x_{k}\right) >2L_{i}$. Since 
\begin{equation*}
u_{i}^{n}\left( x_{k}\right) \leq u_{i}^{n+1}\left( x_{k}\right) \leq
...\leq u_{i}^{n+j}\left( x_{k}\right) \leq ...\leq u_{i}\left( x_{k}\right)
\forall j\text{, }
\end{equation*}%
we have that $u_{i}\left( x_{k}\right) \geq 2L_{i},$ which is a
contradiction. Therefore, $u_{i}$ blows up at the boundary. This solution $%
(u_{1},...,u_{d})$ dominates all other solutions and is therefore commonly
called blow-up/large solution.

\subsubsection{Proof of 3:}

Let $\left( ...,u_{i_{0}}^{n},...,u_{j_{0}}^{n},...\right) \in C^{1}\left(
\Omega \right) $ be the solution of the problem ($P_{\alpha }$) with $\alpha
_{i_{0}}=n$, $n\in N$, and $\alpha _{j_{0}}$ fixed. As in the previous case,
the sequence $u_{i_{0}}^{n}$ is bounded on a compact subset contained in \ $%
\Omega $, implying that there exist functions $u_{i_{0}}$ ($i=1,...,d$)
satisfying $u_{i_{0}}^{n}\rightarrow u_{i_{0}}$ in $C^{1}(K)$ ($i=1,...,d$)
for any compact subset $K\subset \Omega $. Moreover, the arguments used in
the previous cases yield that $u_{i_{0}}$ blows up at the boundary, that is, 
$u_{i_{0}}^{n}=\infty $ on $\partial \Omega $. Related to the sequence ($%
u_{j_{0}}$), we recall that%
\begin{equation*}
\left\{ 
\begin{array}{l}
\Delta _{p}u_{j_{0}}^{n}=F_{u_{j_{0}}}\left(
x,u_{1}^{n},...,u_{j_{0}}^{n},...,u_{d}^{n}\right) \text{ in }\Omega , \\ 
u_{j_{0}}^{n}=\alpha _{j_{0}}\text{ on }\partial \Omega .%
\end{array}%
\right.
\end{equation*}%
Then, by the comparison principle $u_{j_{0}}^{n}\leq \alpha _{j_{0}}$ $%
\forall x\in \overline{\Omega }$ and $n\geq 1$. Passing to the limit as $%
n\rightarrow \infty $, we obtain that $u_{j_{0}}\leq \alpha _{j_{0}}$ for
all $x\in \Omega $.

Claim. Let $x_{0}\in \partial \Omega $ and $\left( x_{k}\right) \subset
\Omega $ be a sequence with $x_{k}\rightarrow x_{0}$. Then $u_{j_{0}}\left(
x_{k}\right) \rightarrow \alpha _{j_{0}}$ as $k\rightarrow \infty $.

Indeed, if the limit does not hold, there exist $\varepsilon >0$ and a
subsequence of $\left( x_{k}\right) $, still denoted by itself, such that%
\begin{equation}
x_{k}\rightarrow x_{0}\text{ and }u_{j_{0}}\left( x_{k}\right) \leq \alpha
_{j_{0}}-\varepsilon \text{ }\forall k\in \mathbb{N}.  \label{3.4}
\end{equation}%
Since $u_{j_{0}}=\alpha _{j_{0}}$ on $\partial \Omega $ and is continuous,
there is some $\delta >0$ such that $u_{j_{0}}\left( x_{k}\right) \geq
\alpha _{j_{0}}-\frac{\varepsilon }{2},$ $\forall x\in \overline{\Omega }%
_{\delta }$. Hence, for $k$ large enough, $x_{k}\in \overline{\Omega }%
_{\delta }$ and $u_{j_{0}}^{1}\geq u_{j_{0}}\left( x_{k}\right) \geq \alpha
_{j_{0}}-\frac{\varepsilon }{2}>\alpha _{j_{0}}-\varepsilon $ which
contradicts (\ref{3.4}). From this claim, we can continuously extend the
function $u_{j_{0}}$ from $\Omega $ to $\overline{\Omega }$ by considering $%
u_{j_{0}}\left( x\right) =\alpha _{j_{0}}$ on $\partial \Omega $, concluding
this way the proof of the Finite and infinite case.

\subsection{Proof of Theorem \protect\ref{entire}}

Firstly, we provide a sub-solution for the problem (\ref{LS}). To do this we
consider the function $w:\mathbb{R}^{N}\rightarrow \left[ 0,\infty \right) $
implicitly defined by%
\begin{equation*}
z\left( x\right) =\int_{w\left( x\right) }^{\infty }\frac{1}{g^{1/\left(
p-1\right) }\left( t\right) }dt,\text{ }x\in \mathbb{R}^{N}.
\end{equation*}%
Note that $w\in C^{1}\left( \mathbb{R}^{N},\left( 0,\infty \right) \right) $%
, $w\left( x\right) \rightarrow +\infty $ as $\left\vert x\right\vert
\rightarrow \infty $ and%
\begin{eqnarray}
\nabla z\left( x\right) &=&-g^{-1/\left( p-1\right) }\left( w\left( x\right)
\right) \nabla w\left( x\right)  \label{p1} \\
\left\vert \nabla w\left( x\right) \right\vert ^{p-2}\nabla w\left( x\right)
&=&-g\left( w\left( x\right) \right) \left\vert \nabla z\left( x\right)
\right\vert ^{p-2}\nabla z\left( x\right) \text{.}  \label{p2}
\end{eqnarray}%
Given $\phi \in C_{0}^{\infty }\left( \mathbb{R}^{N}\right) $, $\phi \geq 0$
we have%
\begin{eqnarray*}
\int_{\mathbb{R}^{N}}\left\vert \nabla w\left( x\right) \right\vert
^{p-2}\nabla w\left( x\right) \nabla \phi dx &=&\int_{\mathbb{R}%
^{N}}-g\left( w\left( x\right) \right) \left\vert \nabla z\left( x\right)
\right\vert ^{p-2}\nabla z\left( x\right) \nabla \phi dx \\
&=&\int_{\mathbb{R}^{N}}\func{div}\left[ g\left( w\left( x\right) \right)
\left\vert \nabla z\left( x\right) \right\vert ^{p-2}\nabla z\left( x\right) %
\right] \phi dx.
\end{eqnarray*}%
Computing the derivatives in the integrand of the expression just above, in
the distribution sense, using (\ref{p1}) we get,%
\begin{eqnarray*}
\int_{\mathbb{R}^{N}}\left\vert \nabla w\left( x\right) \right\vert
^{p-2}\nabla w\left( x\right) \nabla \phi dx &=&\int_{\mathbb{R}^{N}}g\left(
w\left( x\right) \right) \Delta _{p}z\left( x\right) \phi dx \\
&&-\int_{\mathbb{R}^{N}}g^{\prime }\left( w\left( x\right) \right) g^{\frac{1%
}{p-1}}\left( w\left( x\right) \right) \left\vert \nabla z\left( x\right)
\right\vert ^{p}\phi dx.
\end{eqnarray*}%
Using the fact that $g\in \mathcal{F}$ and that $z\left( x\right) $ is an
upper solution of (\ref{qa}) we derive the inequality 
\begin{equation*}
\int_{\mathbb{R}^{N}}\left\vert \nabla w\left( x\right) \right\vert
^{p-2}\nabla w\left( x\right) \nabla \phi dx-\int_{\mathbb{R}^{N}}g\left(
w\left( x\right) \right) \Delta _{p}z\left( x\right) \phi dx\leq 0,
\end{equation*}%
and so%
\begin{equation*}
\int_{\mathbb{R}^{N}}\left\vert \nabla w\left( x\right) \right\vert
^{p-2}\nabla w\left( x\right) \nabla \phi dx+\int_{\mathbb{R}^{N}}g\left(
w\left( x\right) \right) \left( \sum_{i=1}^{d}a_{i}\left( x\right) \right)
\phi dx\leq 0,
\end{equation*}%
which together with (\ref{1.4}) leads to%
\begin{equation*}
-\int_{\mathbb{R}^{N}}\left\vert \nabla w\left( x\right) \right\vert
^{p-2}\nabla w\left( x\right) \nabla \phi dx\geq \int_{\mathbb{R}%
^{N}}a_{i}\left( x\right) F_{u_{i}}\left( x,w\left( x\right) ,...,w\left(
x\right) \right) \phi dx
\end{equation*}%
for all $i=1,...,d$.

In the next, we consider the system%
\begin{equation}
\left\{ 
\begin{array}{l}
\Delta _{p}u_{i}=a_{i}\left( x\right) F_{u_{i}}\left(
x,u_{1},...,u_{i},...,u_{d}\right) \text{ in }B_{n} \\ 
u_{i}=w_{n}\text{ in }\partial B_{n}, \\ 
i=1,...,d,%
\end{array}%
\right.  \label{ball}
\end{equation}%
where $B_{n}$ is the open ball of radius $n$ centered at the origin and $%
w_{n}=\max_{x\in \overline{B}_{n}}w\left( x\right) $. Clearly, $\left(
w,...,w\right) $ and $\left( w_{n},...,w_{n}\right) $ are a sub-solution and
super-solution for (\ref{ball}) respectively. Thus, by Theorem \ref{1.1},
there is a solution $\left( u_{1}^{n},...,u_{d}^{n}\right) \in \left[
W^{1,p}\left( B_{n}\right) \right] ^{d}$ of (\ref{ball}) satisfying $w\left(
x\right) \leq u_{i}^{n}\leq w_{n}$ for all $x\in \overline{B}_{n}$ and $%
i=1,...,d$. For $m\geq 1$ and $n\geq m+1$ consider the family of systems%
\begin{equation*}
\left\{ 
\begin{array}{l}
\Delta _{p}u_{i}^{n}=a_{i}^{a}f_{i}\left( u_{i}^{n}\right) \text{ in }B_{m+1}
\\ 
i=1,...,d,%
\end{array}%
\right.
\end{equation*}%
where $a_{i}^{a}=\min_{x\in \overline{B}_{n}}a_{i}\left( x\right) >0$.
Arguing as in the previous sections, there are monotone decreasing functions 
$\mu _{i}^{a}:\left( 0,\infty \right) \rightarrow \left( 0,\infty \right) $
determined by $f_{i}$ such that 
\begin{equation*}
w\left( x\right) \leq u_{i}^{n}\left( x\right) \leq \mu _{i}^{a}\left(
dist\left( x,\partial B_{m+1}\right) \right) \text{ }\forall x\in B_{m+1}
\end{equation*}%
from which it follows that 
\begin{equation*}
w\left( x\right) \leq u_{i}^{n}\left( x\right) \leq M_{i}^{m}\text{ for all }%
n\in \mathbb{N}\text{, }x\in \overline{B}_{m}\text{, }i=1,...,d\text{ }
\end{equation*}%
and for some positive constants $M_{i}^{m}$. Now using the fact that $%
u_{i}^{n}\in W^{1,p}\left( B_{m}\right) \cap L^{\infty }\left( B_{m}\right) $
if follows from the results of DiBenedetto \cite{dib} and Lieberman \cite%
{lib} that there exist some constants $C_{i}:=C_{i}\left( p,N,\left\vert
u_{i}^{n}\right\vert _{\infty },B_{m}\right) >0$ such that $u_{i}^{n}\in
C^{1,\alpha }\left( B_{m}\right) $ and 
\begin{equation*}
\left\Vert u_{i}^{n}\right\Vert _{C^{1,\alpha }\left( B_{m}\right) }\leq
C_{i}\text{, }i=1,...,d\text{ and }\alpha \in \left( 0,1\right) .\text{ }
\end{equation*}
\ As a consequence, there is $u_{i}\in C^{1}\left( B_{m}\right) $ ($%
i=1,...,d $) such that for some sub-sequence of $u_{i}^{n}$, still denoted
by itself, we get 
\begin{equation*}
u_{i}^{n}\rightarrow u_{i}\text{ }(i=1,...,d)\text{ pointwisely in }B_{m}%
\text{ (}\forall m>1\text{).}
\end{equation*}%
Therefore, $\left( u_{1},...,u_{d}\right) \in C^{1}\left( \mathbb{R}%
^{N}\right) $ and is a solution for the system%
\begin{equation*}
\left\{ 
\begin{array}{l}
\Delta _{p}u_{i}=a_{i}\left( x\right) F_{u_{i}}\left(
x,u_{1},...,u_{i},...,u_{d}\right) \text{ in }\mathbb{R}^{N} \\ 
u_{i}>0\text{ in }\mathbb{R}^{N} \\ 
i=1,...,d%
\end{array}%
\right.
\end{equation*}%
satisfying%
\begin{equation}
w\left( x\right) \leq u_{i}\left( x\right) \text{ for all }x\in \mathbb{R}%
^{N}\text{ and }i=1,...,d.  \label{finish}
\end{equation}%
Letting $\left\vert x\right\vert \rightarrow \infty $ in (\ref{finish}) it
follows that $(u_{1},...,u_{d})$ is a large entire solution for (\ref{LS}).

\section{Remarks}

Assume that $\psi $ belongs to a wide class $\Psi $ of monotone increasing
convex functions. There is an area in probability theory where
boundary-blow-up problem 
\begin{equation*}
\left\{ 
\begin{array}{l}
\Delta u=\psi \left( u\right) \text{ in }\Omega \text{ } \\ 
u=\infty \text{ on }\partial \Omega%
\end{array}%
\right.
\end{equation*}%
arise (see the paper \cite{din1} or directly the book \cite{din2} for
details). The area is known as the theory of superdiffusions, a theory which
provides a mathematical model of a random evolution of a cloud of particles.
Indeed, given any bounded open set $\Omega $ in the N-dimensional Euclidean
space, and any finite measure $\mu $ we may associate with these the exit
measure from $\Omega $ i.e. $\left( X_{\Omega },P_{\mu }\right) $, a random
measure which can be constructed by a passage to the limit from a particles
system. Particles perform independently $\Delta $-diffusions and they
produce, at their death time, a random offspring (cf. \cite{din3}). $P_{\mu
} $ is a probability measure determined by the initial mass distribution $%
\mu $ of the offspring and $X_{\Omega }$ corresponds to the instantaneous
mass distribution of the random evolution cloud. Then procedding in this
way, one can obtain any function $\psi $ from a subclass $\Psi _{0}$ of $%
\Psi $ which contains $u^{\gamma }$ with $1<\gamma \leq 2$. Dynkin \cite%
{din1}, also provided a simple probabilistic representation of the solution
for the class of problems $u^{\gamma }$ ($1<\gamma \leq 2$), in terms of the
so-called exit measure of the associated superprocess. Moreover, the author
say that a probabilistic interpretation is known only for $1<\gamma \leq 2$.

We also remark from the paper of Lasry-Lions \cite{lasry} and Busca-Sirakov 
\cite{busca} that the solutions of the system (\ref{P}) can be viewed as the
value function of a stochastic control process, and the boundary conditions
then means that the process is discouraged to leave the domain by setting an
infinite cost on the boundary. For a more detailed discussion about
practical applications where such problems appear we advise the reader the
introduction of the work \cite{matero}.

\end{document}